# SPECIAL ALGEBRAIC STRUCTURES

**by Florentin Smarandache**
**University of New Mexico**
**200 College Raod**
**Gallup, NM 87301, USA**

**Abstract**.
New notions are introduced in algebra in order to better study the congruences in number theory.
For example, the <special semigroups> make an important such contribution.

**Introduction**.
By <proper subset> of a set A we consider a set P included in A,
and different from A, different from the empty set, and from the unit element
in A - if any.

We rank the algebraic structures using an order relationship:
we say that the algebraic structures S1 << S2 if:
- both are defined on the same set;
- all S1 laws are also S2 laws;
- all axioms of an S1 law are accomplished by the corresponding S2 law;
- S2 laws accomplish strictly more axioms than S1 laws, or S2 has more laws
than S1.

For example: semigroup << monoid << group << ring << field,
or semigroup << commutative semigroup, ring << unitary ring, etc.

We define a **GENERAL SPECIAL STRUCTURE** to be a structure SM on a set A,
different from a structure SN, such that a proper subset of A is an SN
structure, where SM << SN.

1) **The SPECIAL SEMIGROUP** is defined to be a semigroup A, different from
a group, such that a proper subset of A is a group (with respect to the same
induced operation).

For example, if we consider the commutative multiplicative group
  SG = {18^2, 18^3, 18^4, 18^5} (mod 60)
we get the table:

```
   x  |  24 12 36 48
  --- |-------------
   24 |  36 48 24 12
   12 |  48 24 12 36
   36 |  24 12 36 48
   48 |  12 36 48 24
```

Unitary element is 36.

Using the algorithm [Smarandache 1972] we get that
    18^2 is congruent to 18^6 (mod 60).

Now we consider the commutative multiplicative semigroup
    SS = {18^1, 18^2, 18^3, 18^4, 18^5} (mod 60)
and we get the table:

```
    x | 18 | 24 12 36 48
  ----|----|------------
   18 | 24 | 12 36 48 24
  ----|----|------------
   24 | 12 | 36 48 24 12
   12 | 36 | 48 24 12 36
   36 | 48 | 24 12 36 48
   48 | 24 | 12 36 48 24
```

Because SS contains a proper subset SG, which is a group, then SS is a
Special Semigroup.   This is generated by the element 18.  The
powers of 18 form a cyclic sequence: 18,  24,12,36,48,  24,12,36,48, ... .

Similarly are defined:

2) **The SPECIAL MONOID** is defined to be a monoid A, different from a group,
such that a proper subset of A is a group (with respect with the same induced
operation).

3) **The SPECIAL RING** is defined to be a ring A, different from a field, such
that a proper subset of A is a field (with respect with the same induced
operations).

We consider the commutative additive group M={0,18^2,18^3,18^4,18^5}
(mod 60) [using the module 60 residuals of the previous powers of 18],
M={0,12,24,36,48}, unitary additive unit is 0.
(M,+,x) is a field.
While (SR,+,x)={0,6,12,18,24,30,36,42,48,54} (mod 60) is a ring whose

proper subset {0,12,24,36,48} (mod 60) is a field.
Therefore (SR,+,x) (mod 60) is a Special Ring.
This feels very nice.

4) **The SPECIAL SUBRING** is defined to be a Special Ring B which
is a proper subset of a Special Ring A (with respect with the same
induced operations).

5) **The SPECIAL IDEAL** is defined to be an ideal A, different from a field,
such that a proper subset of A is a field (with respect with the same induced
operations).

6) **The SPECIAL SEMILATTICE** is defined to be a lattice A, different from a
lattice, such that a proper subset of A is a lattice (with respect with the
same induced operations).

7) **The SPECIAL FIELD** is defined to be a field (A,+,x), different from a
K-algebra, such that a proper subset of A is a K-algebra (with respect with
the same induced operations, and an external operation).

8) **The SPECIAL R-MODULE** is defined to be an R-MODULE (A,+,x), different from
an S-algebra, such that a proper subset of A is an S-algebra (with respect
with the same induced operations, and another "x" operation internal on A),
where R is a commutative unitary ring and S is its proper subset field.

9) **The SPECIAL K-VECTORIAL SPACE** is defined to be a K-vectorial
space (A,+,.), different from a K-algebra, such that a proper subset of A is
a K-algebra (with respect with the same induced operations, and another "x"
operation internal on A), where K is a commutative field.

Craiova,  1973

and Applied Sciences>, Delhi, India, Vol. 17E, No. 1, 119-121, 1998.